\documentclass[10pt, a4paper]{amsart}
\usepackage{amssymb,amsmath,amsthm,latexsym,amscd,amsfonts,mathrsfs}
\usepackage{graphicx}
\usepackage[utf8]{inputenc}

\newcommand{\bd}{{\mathbb{D}}}

\newcommand{\br}{{\mathbb{R}}}

\newcommand{\bc}{{\mathbb{C}}}

\newcommand{\s}{\sigma}

\newcommand{\p}{\varphi}

\newcommand{\oo}{\Omega}

\newcommand{\gga}{\Gamma}
\newcommand{\ep}{\varepsilon}
\newcommand{\z}{\zeta}

\newcommand{\nt}{\noindent}

\newcommand{\bsl}{\backslash}

\newcommand{\lp}{\left(}
\newcommand{\rp}{\right)}

\DeclareMathOperator{\dist}{{\it d}}

\newcommand{\Sc}{\mathcal{S}}

\allowdisplaybreaks
\numberwithin{equation}{section}

\newtheorem{theorem}{Theorem}[section]
\newtheorem*{theorem-a}{Theorem A}
\newtheorem*{theorem-b}{Theorem B}

\newtheorem{remark}[theorem]{Remark}

\begin{document}

\title[Fredholm alternative]
{A remark on analytic Fredholm alternative}

\author{L. Golinskii}
\address{Mathematics Division, ILTPE, 47 Science ave.,  61103 Kharkov, Ukraine}
\email{golinskii@ilt.kharkov.ua}

\author{S. Kupin}
\address{IMB, Universit\'e de Bordeaux, 351 ave. de la Lib\'eration, 33405 Talence Cedex, France}
\email{skupin@math.u-bordeaux1.fr}

\subjclass{Primary: 47A56; Secondary: 47B10}

\begin{abstract}
We apply a recent result of Borichev--Golinskii--Kupin on the Blasch\-ke-type conditions for zeros of analytic
functions on the complex plane with a cut along the positive semi-axis to the problem of the eigenvalues
distribution of the Fredholm-type analytic operator-valued functions.
\end{abstract}

\maketitle

\vspace{0.5cm}
\section*{Introduction and main results}\label{s0}

The goal of this note is to refine partially (for a certain range of parameters) a
recent result of R. Frank \cite[Theorem 3.1]{fra15} on some quantitative aspects of the analytic Fredholm alternative.
Precisely, the problem concerns the distribution of eigenvalues of finite type of an
operator-valued function $W(\cdot)=I+T(\cdot)$, analytic on a domain $\oo$ of the complex plane.
We always assume that $T\in\Sc_\infty$, the set of compact operators on the Hilbert space.
A  number $\lambda_0\in\oo$ is called an {\it eigenvalue of finite type of} $W$ if $\ker W(\lambda_0)\not=\{0\}$,
(i.e., $-1$ is an eigenvalue of $T(\lambda_0)$), if $W(\lambda_0)$ is Fredholm (that is, both
$\dim\ker W(\lambda_0)$ and ${\rm codim\ ran}\, W(\lambda_0)$ are finite), and if $W$ is invertible in some
punctured neighborhood of $\lambda_0$. The function $W$ admits the following
expansion at any eigenvalue of finite type, see \cite[Theorem XI.8.1]{ggk},
$$ W(\lambda)=E(\lambda)(P_0+(\lambda-\lambda_0)^{k_1}\,P_1+\ldots+(\lambda-\lambda_0)^{k_l}\,P_l)G(\lambda), $$
where $P_1,\ldots,P_l$ are mutually disjoint projections of rank one, $P_0=I-P_1-\ldots-P_l$,
$k_1\le\ldots\le k_l$ are positive integers, and $E$, $G$ are analytic operator-valued functions,
defined and invertible in some neighborhood of $\lambda_0$. The number
$$ \nu(\lambda_0,W):=k_1+\ldots+k_l $$
is usually referred to as an {\it algebraic multiplicity of the eigenvalue} $\lambda_0$.

The following result, Theorem 3.1, is a cornerstone of the paper \cite{fra15}.
By $\{\lambda_j\}$ we always denote the eigenvalues of $W=I+T$ of finite type, repeated accordingly to their algebraic multiplicity.

\begin{theorem-a}\label{th-a}
Let $T(\cdot)$ be an analytic operator-valued function on the domain $\oo=\bc\backslash\br_+$, so that $T\in\Sc_p$,
$p\ge 1$, the set of the Schatten--von Neumann operators of order $p$. Assume that for all $\lambda\in\bc\backslash\br_+$
\begin{equation}\label{eq1}
\|T(\lambda)\|_p\le \frac{M}{\dist^\rho(\lambda,\br_+)|\lambda|^\s}\,, \qquad  \rho>0, \ \ \s\in\br, \ \ \rho+\s>0,
\end{equation}
$\dist(\lambda,\br_+)$ is the Euclidean distance from $\lambda$ to the positive semi-axis.
Then for all $\ep,\,\ep'>0$ and $\nu\ge1$
\begin{equation}\label{freq2}
\sum_{|\lambda_j|\le M^{1/(\rho+\s)}}  \dist^{p\rho+1+\ep}(\lambda_j,\br_+)\, |\lambda_j|^{\frac{q-p\rho-1-\ep}2}\le
C M^{\frac{q+p\rho+1+\ep}{2(\rho+\s)}},
\end{equation}
where $q:=(p\rho+2p\s-1+\ep)_+$, and
\begin{equation}\label{freq1}
\sum_{|\lambda_j|\ge\nu M^{1/(\rho+\s)}}  \dist^{p\rho+1+\ep}(\lambda_j,\br_+)\, |\lambda_j|^{\rho+\s-p\rho-1-\ep-\ep'}\le
\frac{C}{\nu^{\ep'}}\, M^{\frac{\rho+\s-\ep'}{\rho+\s}}.
\end{equation}
Here $C$ is a generic positive constant which depends on $p, \rho, \s, \ep, \ep'$.
\end{theorem-a}

The similar results for $\rho=0$ are also available.

The proof of this result is based on the identification of the eigenvalues of finite type of $W$ with
the zeros of certain {\it scalar} analytic functions, known as the {\it regularized determinants}
$$ f(\lambda):={\rm det}_{p}(I+T(\lambda)), $$
see \cite{gokr-nsa, Si05} for their definition and basic properties. The point is that the set of eigenvalues
of finite type of $W$ agrees with the zero set of $f$, and moreover, $\nu(\lambda_0,W)=\mu_f(\lambda_0)$, the multiplicity of zero of
$f$ at $\lambda_0$ (see \cite[Lemma 3.2]{fra15} for the rigorous proof). Thereby, the problem is reduced to the study of the zero
distributions of certain analytic functions, the latter being a classical topic of complex analysis going back to Jensen \cite{jen} and
Blaschke \cite{bla}.

A key ingredient of the proof in \cite{fra15} is a result of \cite[Theorem 0.2]{bgk1} on the Blaschke-type conditions for zeros of
analytic functions in the unit disk which can grow at the direction of certain (finite) subsets of the unit circle. In a recent manuscript
\cite{bgk2} some new such conditions on zeros of analytic functions in the unit disk and on some other domains, including the complex
plane with a cut along the positive semi-axis, are suggested. Here is a particular case of \cite[Theorem 4.5]{bgk2} which seems relevant.
We use a convenient shortening
\begin{equation*}
\{u\}_{c,\ep}:=(u_- -1+\ep)_+ -\min(c,u_+), \qquad c\ge0, \ \ \ep>0, \ \ u=u_+-u_-\in\br.
\end{equation*}

\begin{theorem-b}
Let $h$ be an analytic function on $\oo=\bc\backslash\br_+$, $|h(-1)|=1$, subject to the growth condition
\begin{equation*}
\log|h(\lambda)|\le \frac{K}{|\lambda|^r}\,\frac{(1+|\lambda|)^{b}}{\dist^a(\lambda,\br_+)}, \quad \lambda \in\bc\bsl\br_+, \quad a,b\ge0, \ \ r\in\br.
\end{equation*}
Let $Z(h)$ be its zero set, counting multiplicities (the divisor of $h$). Denote
$$ s:=3a-2b+2r. $$
Then for each $\ep>0$ there is a positive number $C$ which depends on all parameters
involved such that the following inequality holds
\begin{equation}\label{btc08}
\sum_{z\in Z(h)} \dist^{a+1+\ep}(z,\br_+)\,\frac{|z|^{s_1}}{(1+|z|)^{s_2}}\le C\cdot K,
\end{equation}
where the parameters $s_1$, $s_2$ are defined by the relations
\begin{equation*}
s_1 := \frac{\{-2r-a\}_{a,\ep}-a-1-\ep}2\,, \qquad s_2 := a+1+\ep + \frac{\{-2r-a\}_{a,\ep}+\{s\}_{a,\ep}}2\,.
\end{equation*}
\end{theorem-b}

We are aimed at proving the results, which refine Theorem A for a certain range of parameters, by using Theorem B.

\begin{theorem}\label{t1}
Let $T(\cdot)$ be an analytic operator-valued function on the domain $\oo=\bc\backslash\br_+$,
which satisfies the hypothesis of Theorem A. Assume that
\begin{equation}\label{range1}
0<\rho+\s\le\frac{\rho}2\,.
\end{equation}
Then for all $0<\ep<1$
\begin{equation}\label{gkeq1}
\sum_{|\lambda_j|\le M^{1/(\rho+\s)}}  \dist^{p\rho+1+\ep}(\lambda_j,\br_+)\, |\lambda_j|^{p\s-\frac{1+\ep}2}\le
C M^{p+\frac{1+\ep}{2(\rho+\s)}}.
\end{equation}
\end{theorem}

\smallskip

Note that under assumption \eqref{range1}
$$ p\s-\frac{1+\ep}2\le -\frac{p\rho+1+\ep}2<0, $$
so for $|\z|\le1$
$$ |\z|^{p\s-\frac{1+\ep}2}\ge |\z|^{-\frac{p\rho+1+\ep}2}\ge |\z|^{\frac{q-p\rho-1-\ep}2}\,, $$
that is, \eqref{gkeq1} is stronger than \eqref{freq2} with regard to eigenvalues tending to zero. Theorem B
gives the same results, \eqref{freq2} and \eqref{freq1}, as in Theorem A, for the rest of the values of $\rho$ and $\s$, and the
eigenvalues  tending to infinity.

\smallskip

The case
\begin{equation}\label{range11}
\rho>0, \qquad \rho+\s<0,
\end{equation}
is not treated in \cite{fra15}.

\begin{theorem}\label{t2}
Under conditions $\eqref{range11}$ assume that for all $\lambda\in\bc\backslash\br_+$
\begin{equation}\label{eq11}
\|T(\lambda)\|_p\le \frac{M}{\dist^\rho(\lambda,\br_+)|\lambda|^\s}\,.
\end{equation}
Then for  $-\rho/2\le\rho+\s<0$ and all $\ep>0$
\begin{equation}\label{gkeq2}
\sum_{|\lambda_j|\ge M^{1/(\rho+\s)}}  \dist^{p\rho+1+\ep}(\lambda_j,\br_+)\, |\lambda_j|^{p\s-\frac32\,(1+\ep)}\le
C\, M^{p-\frac{1+\ep}{2(\rho+\s)}}\,,
\end{equation}
and for $\rho+\s<-\rho/2$ and all $\ep>0$
\begin{equation}\label{gkeq4}
\sum_{|\lambda_j|\ge M^{1/(\rho+\s)}}  \dist^{p\rho+1+\ep}(\lambda_j,\br_+)\, |\lambda_j|^{-\frac{l+3(p\rho+1+\ep)}2}\le
C\, M^{-\frac{l+p\rho+1+\ep}{2(\rho+\s)}}\,,
\end{equation}
where $l:=(-3p\rho-2p\s-1+\ep)_+$. Moreover, under conditions $\eqref{range11}$, for all $\ep, \ep'>0$ and $0<\mu\le1$
\begin{equation}\label{gkeq3}
\sum_{|\lambda_j|\le\mu M^{1/(\rho+\s)}}  \dist^{p\rho+1+\ep}(\lambda_j,\br_+)\, |\lambda_j|^{\rho+\s-p\rho-1-\ep+\ep'}\le
C\,\mu^{\ep'}\,M^{\frac{\rho+\s+\ep'}{\rho+\s}}\,.
\end{equation}
\end{theorem}

\section{Proof of main results}\label{s1}

\smallskip
\nt
{\it Proof of Theorem \ref{t1}.}

\nt
We follow the line of reasoning from \cite{fra15}. The scaling $T_1(\lambda):=T(M^{1/(\rho+\s)}\,\lambda)$ looks reasonable, so
$$ \|T_1(\lambda)\|_p\le \frac{1}{\dist^\rho(\lambda,\br_+)\,|\lambda|^\s}\,, $$
and, by \cite[Theorem 9.2, (b)]{Si05}, we have for the determinant $f_1=\det_{p}(I+T_1)$
\begin{equation}\label{eq2}
\log|f_1(\lambda)|\le \frac{\gga_p}{\dist^{p\rho}(\lambda,\br_+)\,|\lambda|^{p\s}}\,, \qquad \lambda\in\bc\backslash\br_+.
\end{equation}

To apply Theorem B we have to ensure the normalization condition. Note that the function $T_1$ tends to zero along the
left semi-axis as long as $\rho+\s>0$, so the inequality
(see \cite[Theorem 9.2, (c)]{Si05})
\begin{equation}\label{belbound}
|f_1(\lambda)-1|\le \p(\|T_1(\lambda)\|_p), \quad \p(t):=t\exp\bigl(\gga_p(t+1)^p\bigr), \quad t\ge0,
\end{equation}
holds with a suitable constant $\gga_p$ which depends only on $p$, and provides a lower bound for $f_1$ whenever the right
side is small enough. We have for $t\ge1$ and $\lambda=-t\in\br_-$
$$ |f_1(-t)-1|\le \frac{C_1}{t^{\rho+\s}}\,, $$
(in the sequel $C_k$ stand for generic positive constants depending on the parameters involved). If $t\ge (2C_1)^{1/(\rho+\s)}=C_2$, then
$|f_1(-t)|\ge 1/2$, and so
\begin{equation}\label{eq3}
\log|f_1(-t)|\ge -2(1-|f_1(-t)|)\ge -\frac{2C_1}{t^{\rho+\s}}\,.
\end{equation}

Next, put
\begin{equation}\label{norm}
h(\lambda):=\frac{f_1(t\l)}{f_1(-t)}\,, \qquad h(-1)=1.
\end{equation}
It follows from \eqref{eq2} and \eqref{eq3} that for $t\ge C_2$
\begin{equation*}
\begin{split}
\log|h(\lambda)| &= \log|f_1(t\l)|-\log|f_1(-t)|\le
\frac{\gga_p}{t^{\rho+\s}}\,\frac{1}{\dist^{p\rho}(\lambda,\br_+)\,|\lambda|^{p\s}}+\frac{2C_1}{t^{\rho+\s}} \\
&\le \frac{C_3}{t^{\rho+\s}}\,\lp\frac{1}{\dist^{p\rho}(\lambda,\br_+)\,|\lambda|^{p\s}}+1\rp
\le  \frac{C_3}{t^{\rho+\s}}\,\frac{(1+|\lambda|)^{p(\rho+\s)}}{\dist^{p\rho}(\lambda,\br_+)\,|\lambda|^{p\s}}\,.
\end{split}
\end{equation*}
Theorem B applies now with
$$ a=p\rho, \quad r=p\s, \quad b=p(\rho+\s), \quad K=\frac{C_3}{t^{\rho+\s}}\,, $$
and $s=a$, $\{s\}_{a,\ep}=-a$. In view of \eqref{range1} one has $2r+a=p(\rho+2\s)\le0$, so
$$ \{-2r-a\}_{a,\ep}=-\min(a, -2r-a)=2r+a=p\rho+2p\s, $$
(recall that, by the assumption, $a>-2r-a$). Hence
$$ s_1=\frac{2p\s-1-\ep}2\,, \qquad s_2=p\rho+p\s+1+\ep, $$
and \eqref{btc08} implies
\begin{equation*}
\sum_{z\in Z(h)} \dist^{p\rho+1+\ep}(z,\br_+)\,\frac{|z|^{\frac{2p\s-1-\ep}2}}{(1+|z|)^{p\rho+p\s+1+\ep}}
\le \frac{C_4}{t^{\rho+\s}}\,,
\end{equation*}
or
\begin{equation*}
t^{\frac{1+\ep}2}\,\sum_{\z\in Z(f_1)} \dist^{p\rho+1+\ep}(\z,\br_+)\,\frac{|\z|^{\frac{2p\s-1-\ep}2}}{(t+|\z|)^{p\rho+p\s+1+\ep}}
\le \frac{C_4}{t^{\rho+\s}}\,.
\end{equation*}

For $|\z|\le1$ we fix $t$, say, $t=C_2$, and since $t+|\z|\le C_2+1$, we come to
\begin{equation*}
\sum_{\z\in Z(f_1)\cap\bar\bd} \dist^{p\rho+1+\ep}(\z,\br_+)\,|\z|^{p\s-\frac{1+\ep}2}\le C_5,
\end{equation*}
which, after scaling, is \eqref{gkeq1}. The proof is complete. \hfill $\Box$

\medskip
\nt
{\it Proof of Theorem \ref{t2}.}

\nt
The idea is much the same with the only technical differences. In the above notation relation \eqref{eq2} still holds, and the function
$T_1$ tends to zero as $t\to 0-$ whenever $\rho+\s<0$. So
\begin{equation}\label{eq4}
\log|f_1(-t)|\ge -2(1-|f_1(-t)|)\ge -\frac{2C_1}{t^{\rho+\s}}=-2C_1t^{|\rho+\s|}\,, \qquad 0<t\le C_2.
\end{equation}
For the function $h$ \eqref{norm} we now have
\begin{equation*}
\log|h(\lambda)|\le C_3 t^{|\rho+\s|}\,\lp\frac{1}{\dist^{p\rho}(\lambda,\br_+)\,|\lambda|^{p\s}}+1\rp,
\end{equation*}
and as
$$ \frac{1}{\dist^{p\rho}(\lambda,\br_+)\,|\lambda|^{p\s}}+1\le \frac{|\lambda|^{p|\s|}+|\lambda|^{p\rho}}{\dist^{p\rho}(\lambda,\br_+)}\le |\lambda|^{p\rho}\,
\frac{(1+|\lambda|)^{p|\rho+\s|}}{\dist^{p\rho}(\lambda,\br_+)}\,, $$
we come to the bound
\begin{equation}\label{eq5}
\log|h(\lambda)|\le C_3t^{|\rho+\s|}\,|\lambda|^{p\rho}\,\frac{(1+|\lambda|)^{p|\rho+\s|}}{\dist^{p\rho}(\lambda,\br_+)}\,.
\end{equation}

Theorem B applies with
$$ a=p\rho, \quad r=-a=-p\rho, \quad b=-p(\rho+\s), \quad K=\frac{C_3}{t^{\rho+\s}}\,, $$
and $-2r-a=a>0$, so
$$ \{-2r-a\}_{a,\ep}=-a=-p\rho, \qquad s_1=-p\rho-\frac{1+\ep}2\,. $$

The sign of $s=3a-2b+2r=p(3\rho+2\s)$ (which can be either positive or negative) affects the computation of $\{s\}_{a,\ep}$, so we will
differ two situations. In the case $-\rho/2\le\rho+\s<0$ we have
\begin{equation}\label{eq6}
s>0, \qquad \{s\}_{a,\ep}=-\min(a,s_+)=-s,
\end{equation}
since, by \eqref{range11},  $s_+=s=p(3\rho+2\s)<p\rho=a$. So $s_2=-p(\rho+\s)+1+\ep$, and \eqref{btc08} leads to
\begin{equation}\label{eq7}
t^{p|\rho+\s|+\frac{1+\ep}2}\,\sum_{\z\in Z(f_1)} \dist^{p\rho+1+\ep}(\z,\br_+)\,
\frac{|\z|^{-p\rho-\frac{1+\ep}2}}{(t+|\z|)^{p|\rho+\s|+1+\ep}} \le \frac{C_4}{t^{\rho+\s}}\,, \ \ 0<t\le C_2.
\end{equation}

A simple bound $(C_2+|\z|)^{-1}\ge C_5\,|\z|^{-1}$ for $|\z|\ge1$ and fixed $t=C_2$ gives
\begin{equation*}
\sum_{\z\in Z(f_1)\cap\bd_-} \dist^{p\rho+1+\ep}(\z,\br_+)\,|\z|^{p\s-\frac32\,(1+\ep)}\le C_6, \qquad \bd_-:=\{|\z|\ge1\},
\end{equation*}
which, after scaling, is \eqref{gkeq2}.

If $|\z|\le \mu \le1$, we multiply \eqref{eq7} through by $t^{\rho+\s-1+\ep'}$ and integrate it  termwise with respect to $t$ from $0$ to
$\mu C_2$ (the idea comes from \cite{DeHaKa})
\begin{equation*}
\begin{split}
\int_0^{\mu C_2} \frac{t^{(p-1)|\rho+\s|+\frac{1+\ep}2-1+\ep'}}{(t+|\z|)^{p|\rho+\s|+1+\ep}}\,dt &=|\z|^{\rho+\s-\frac{1+\ep}2+\ep'}\,
\int_0^{\mu C_2/|\z|}\frac{x^{(p-1)|\rho+\s|+\frac{1+\ep}2-1+\ep'}}{(1+x)^{p|\rho+\s|+1+\ep}}\,dx \\
&\ge C_7\,|\z|^{\rho+\s-\frac{1+\ep}2+\ep'}\,,
\end{split}
\end{equation*}
to obtain
\begin{equation*}
\sum_{\z\in Z(f_1)\cap\bd_\mu} \dist^{p\rho+1+\ep}(\z,\br_+)\,|\z|^{\rho+\s-p\rho-1-\ep+\ep'}\le C_8\, (\mu C_2)^{\ep'},
\qquad \bd_\mu:=\{|\z|\le\mu\},
\end{equation*}
which, after scaling, gives \eqref{gkeq3}.

In the case $\rho+\s<-\rho/2$ the proof is the same with $s\le0$ and
$$ \{s\}_{a,\ep}=(-3p\rho-2p\s-1+\ep)_+=l, \qquad s_2=\frac{p\rho+l}2+1+\ep. $$
\hfill $\Box$

\smallskip

\begin{remark}
The case $\eqref{range11}$ can be reduced to the one considered in Theorem A by means of the transformation (the change of variables)
$\lambda\to 1/\lambda$ and the general formula
$$ \dist(1/\lambda, \br_+)=\frac{\dist(\lambda, \br_+)}{|\lambda|^2}\,, \qquad \lambda\in\bc\backslash\br_+. $$
\end{remark}

\smallskip

\begin{remark}
The general form of \cite[Theorem 4.5]{bgk2} allows a finite number of singularities on $\br_+$.
So we can obtain the similar results on eigenvalues of finite type for analytic operator-valued
functions $W=I+T$ on $\bc\backslash\br_+$ subject to the bound
$$ \|T(\lambda)\|_p\le M\,\frac{(1+|\lambda|)^\tau}{\dist^\rho(\lambda,\br_+)|\lambda|^\s}\,\frac{\prod_{j=1}^n |\lambda-t_j|^{c_j}}
{\prod_{k=1}^m |\lambda-t'_k|^{d_k}}\,, \quad \rho,\tau,c_j,d_k\ge0, \ \ \s\in\br, $$
where $\{t_j\}$ and $\{t_k'\}$ are two disjoint finite sets of distinct positive numbers.
\end{remark}

\end{document}